\documentclass[14pt]{article}
\begin {document}
\newtheorem{theorem}{Theorem}[section]
\newtheorem{Prop}{Proposition}

\newtheorem{lemma}[theorem]{Lemma}
\newtheorem{fact}{Example}

\newenvironment{proof}{
\par
\noindent {\bf Proof.}\rm}

\title{Simple decompositions of simple special  Jordan superalgebras}
\author{
M.V.Tvalavadze\thanks{\it Email address: marina@math.mun.ca}, \thanks{supported by NSERC Grant 227060-00}\\
 \small Department of Mathematics and Statistics\\
\small Memorial University of Newfoundland\\
\small St. John's, NL, CANADA. }
\date{}
\maketitle

\sloppy

\noindent-------------------------------------------------------------------------------------------------------\\
{\bf Abstract}

 We classify decompositions of  simple
special finite-dimensional Jordan superalgebras over an
algebraically closed field of characteristic zero  into the sum of
two proper simple subsuperalgebras.\\
\noindent-------------------------------------------------------------------------------------------------------

\section{Introduction}

In our previous work with T.Tvalavadze \cite{tm} we considered
special simple finite-dimensional Jordan algebras decomposable as
the sum of two proper simple subalgebras. The main result in
\cite{tm} is the following.
\par\medskip

{\bf Theorem}. {\it Let $\cal J$ be a finite-dimensional special
simple Jordan algebra over an algebraically closed field $F$ of
 characteristic not two. The only possible decompositions of $\cal J$
as the sum of two simple subalgebras ${\cal J}_1$ and ${\cal J}_2$
are the following:
\par\medskip

{\bf 1}. ${\cal J}\cong F\oplus  V$ and ${\cal J}_1\cong F\oplus
V_1$, ${\cal J}_2\cong F\oplus V_2$, where $V$, $V_1$, $V_2$ are
vector

spaces.

{\bf 2}. Either ${\cal J}\cong  H({\cal R}_3)$ and  ${\cal
J}_1\cong H(F_3)$, ${\cal J}_2\cong F\oplus V$, or ${\cal J}\cong
H( {\cal R}_n)$,

$n\ge 3$, ${\cal J}_1\cong H(F_n)$ and ${\cal J}_2$ is isomorphic
to one of the following

algebras: $H(F_{n-1})$, $H(F_n)$ or $H({\cal R}_{n-1})$.

{\bf 3}. ${\cal J}\cong H({\cal Q}_n)$ and ${\cal J}_1$, ${\cal
J}_2\cong H({\cal R}_n)$.}
\par\medskip
Actually the problem of simple decompositions of simple algebras
first arises in the paper of Onishchik (see \cite{on}) in which he
classified all possible types of simple decompositions of simple
complex and real Lie algebras. Later for associative algebras over
an arbitrary field $F$ the same problem was formulated and then
solved by Bahturin and Kegel in \cite{BK}. According to \cite{BK},
no full-matrix algebra can be written as the sum of two
full-matrix subalgebras. Note that if $F$ is algebraically closed
with zero characteristic, then this follows from \cite{on}.

To begin with we briefly remind the classification of simple
 Jordan superalgebras obtained by Kac (see \cite{Kac}) over an
algebraically closed field $F$ with zero characteristic. If $\cal
J$ is a simple special finite-dimensional Jordan superalgebra over
algebraically closed field $F$ with zero characteristic, then
$\cal J$ is isomorphic to  one of the following superalgebras:
\par\medskip

{\bf (1)} $M_{n,m}(F)^{(+)}$, the set of all matrices of order
$n+m$ with respect to the natural $Z_2$-gradation under the Jordan
supermultiplication;

{\bf (2)} $osp(n,m)$, the set of all matrices of order $n+2m$
symmetric with respect to the orthosymplectic superinvolution. The
superalgebra consists of matrices $\left\{\left(\begin{array}{cc}
                                         A& B\\
                                         C&D
                                         \end{array}\right)\right\}$
where $A^t=A$, $D$ is symplectic, $B$, $C$ are skew-symmetric;

{\bf (3)} $P(n)=\left\{ \left( \begin{array}{cc}
                           A& B\\
                           C& A^t
                           \end{array}\right), B^t=B, C^t=-C\in
                           M_n(F)\right\}$;

{\bf (4)} $Q(n)=\left\{ \left( \begin{array}{cc}
                           A& B\\
                           B& A
                           \end{array}\right)\right\}$
where $A$ and $B$ are any square matrices of order $n$.

{\bf (5)} Let $V=V_0+V_1$ be a $Z_2$-graded vector space with a
non-singular symmetric bilinear superform   $f(\,,\,):V\times V\to
F$. Consider the direct sum of  $F$ and $V$, ${\cal J}=F\oplus V$,
and determine multiplication according to
$$(\alpha+v)(\beta +w)=(\alpha\beta+f(v,w))+(\alpha w+\beta v).$$
Then $\cal J$ becomes a  Jordan superalgebra of the type $J(V,f)$
with respect to the following  $Z_2$-gradation: ${\cal
J}_0=F+V_0$, ${\cal J}_1=V_1$.

{\bf (6)} The 3-dimensional Kaplansky superalgebra  $K_3$,
$(K_3)_0=Fe$, $(K_3)_1=Fx+Fy$, with the multiplication $e^2=e$,
$ex=\frac{1}{2}x$, $ey=\frac{1}{2}y$, $[x,y]=e$.

{\bf (7)} The 1-parametric family of 4-dimensional superalgebras
$D_t$, $D_t=(D_t)_0+(D_t)_1$, where $(D_t)_0=Fe_1+Fe_2$,
$(D_t)_1=Fx+Fy$, where $e_i^2=e_i$, $e_1e_2=0$,
$e_ix=\frac{1}{2}x$, $e_iy=\frac{1}{2}y$, $[x,y]=e_1+te_2$,
$i=1,2$. A superalgebra $D_t$ is simple only if  $t\ne 0$. If
$t=-1$, then  $D_{-1}$ isomorphic $M_{1,1}(F)$.
\par\medskip

Next we cite some important Lemmas and Theorems from \cite{tm}
which will be repeatedly used later.
\par\medskip
\begin{lemma}
Let a Jordan algebra $\cal J$ of the type $H({\cal D}'_m)$ be a
proper subalgebra of $H({\cal D}_n)$ such that the identity of
$H({\cal D}_n)$ is an element  of this subalgebra. If either

\noindent 1. ${\cal D}'=F$ and ${\cal D}=F$, or

\noindent 2. ${\cal D}'={\cal R}, {\cal Q}$, and ${\cal D}={\cal
R} $,

\noindent  then $m\le\frac{n}{2}$
\end{lemma}
\begin{lemma}

Let $V$ be a  vector space with a non-singular symmetric bilinear
form $f$, and $v_0$  a fixed non-trivial vector in $V$. Let $\cal
S$ be the set of  all linear  operators which are symmetric with
respect to  $f$. Then, ${\cal S}v_0=V$.

\end{lemma}
\begin{theorem}

Let $\cal J$ be a simple Jordan algebra of the type $H(F_n)$,  and
$\cal A$, $\cal B$  proper simple  Jordan subalgebras of $\cal J$.
Then ${\cal J}\ne {\cal A}+{\cal B}$.

\end{theorem}

\begin{lemma}
Let $W$ be the natural $2m$-dimensional module for $H({\cal
Q}_m)$, and $v$ be an arbitrary non-zero vector in $W$. Then,
$\dim H({\cal Q}_m)v=2m-1$.
\end{lemma}

Throughout the paper the basic field $F$ is algebraically closed
with characteristic zero.

\section{ Decompositions of  superalgebras of the type
$M_{n,m}(F)^{(+)}$ }

Our main goal is to prove the following.

\begin{theorem}
Let $\cal A$ be a superalgebra of the type $M_{n,m}(F)^{(+)}$
where $n,m>0$. If both $n,m$ are odd, then $\cal A$ has no
decompositions into the sum of two proper nontrivial simple
subsuperalgebras. If one of the indices, for example $m$, is an
even number and the other is odd, then the only possible simple
decomposition is the following: ${\cal A}={\cal B}+{\cal C}$ where
$\cal B$ and $\cal C$ have types $osp(n,\frac{m}{2})$ and
$M_{n-1,m}(F)^{(+)}$, respectively. If both indices are even, then
$\cal A$ admits two  types of decompositions of the following
forms:

1. ${\cal A}={\cal B}_1+{\cal C}_1$ where ${\cal B}_1$ and ${\cal
C}_1$ have types $osp(n,\frac{m}{2})$ and $M_{n-1,m}(F)^{(+)}$,

2. ${\cal A}={\cal B}_2+{\cal C}_2$ where ${\cal B}_2$ and ${\cal
C}_2$ have types $osp(m,\frac{n}{2})$ and $M_{m-1,n}(F)^{(+)}$.
\end{theorem}

Before the discussion of various properties of $M_{n,m}(F)^{(+)}$
we recall a definition of the universal associative enveloping
superalgebra of a Jordan superalgebra which will be frequently
used later.

An associative specialization $u:{\cal J}\to U({\cal J})$ where
$U({\cal J})$ is an associative superalgebra is said to be
universal if $U({\cal J})$ is generated by $u({\cal J})$, and for
any other specialization $\varphi:{\cal J}\to {\cal A}$ where
$\cal A$ is an associative superalgebra there exists a
homomorphism $\psi:U({\cal J})\to {\cal A}$ such that
$\varphi=\psi u$. Then $U({\cal J})$ is called a universal
associative enveloping superalgebra of $\cal J$. It is worth
noting that an associative superalgebra  can be considered as an
associative algebra. The following Theorem by C.Martinez and
E.Zelmanov \cite{MZ} plays a key role in the later discussion.

\begin{theorem} Let $U({\cal J})$ denote a universal associative
enveloping superalgebra for a Jordan superalgebra $\cal J$. Then
$U(M_{k,l}^{(+)})\cong M_{k,l}(F)\oplus M_{k,l}(F)$ where
$(k,l)\ne (1,1)$; $U(Q(k))\cong Q(k)\oplus Q(k)$, $k\ge 2$;
$U(osp(m,n))\cong M_{m, 2n}(F)$, $(m,n)\ne (1,1)$; $U(P(n))\cong
M_{n, n}(F)$, $n\ge 3$.
\end{theorem}
\par\medskip
 {\bf Remark 1}(see \cite{MZ}) {\it In the case where ${\cal J}\cong M_{1,1}(F)^{(+)}, P(2), osp(1,1),
 K_3$ or $D_t$ the universal enveloping superalgebras have more
 complicated structure. Indeed, the universal associative enveloping
 superalgebras of the above Jordan superalgebras are no more
 finite-dimensional. Also we note that if the characteristic of the
 basic field $F$ equals zero, then $K_3$ has no non-zero
 finite-dimensional associative specializations.}

The following Theorem by Martinez and Zelmanov (see \cite{MZ})
describes all irreducible one-sided bimodules of $D(t)$ where
$t\ne -1,0,1$.

\begin{theorem} Let F be an algebraically closed field with zero
characteristic. If $t=-\frac{m}{m+1}$, $m\ge 1$, then $D(t)$ has
two irreducible finite-dimensional one sided bimodules $V_1(t)$
and $V_1(t)^{op}$.

If $t=-\frac{m+1}{m}$, $m\ge 1$, then $D(t)$ has two irreducible
finite-dimensional one sided bimodules $V_2(t)$ and $V_2(t)^{op}$.

If $t$ cannot be represented as $-\frac{m}{m+1}$ and
$-\frac{m+1}{m}$ where $m$ is a positive integer, then $D(t)$ does
not have non-zero finite-dimensional specializations.

\end{theorem}
\par\medskip
{\bf Remark 2} {\it If $\mbox{char}\,F=p>2$, then for an arbitrary
$t$ the superalgebra $D(t)$ can be embedded in a
finite-dimensional associative superalgebra.}

\par\medskip
{\bf Remark 3}  {\it If $t=-\frac{m}{m+1}$, then
$\dim\,V_1(t)_0=m$, $\dim\,V_1(t)_1=m+1$. If $t=-\frac{m+1}{m}$,
then $\dim\,V_2(t)_0=m+1$, $\dim\,V_2(t)_1=m$.}
\par\medskip

Now we look at the case when ${\cal J}\cong J(V,f)$. Let
$V=V_0+V_1$ be a $Z_2$-graded vector space, $\dim\,V_0=m$,
$\dim\,V_1=2n$. Let $f(\,,\,):V\times V\to F$ be a supersymmetric
bilinear form on $V$. The universal associative enveloping algebra
of the Jordan algebra $F+V_0$ is the Clifford algebra
$C(V_0,f)=\langle 1,e_1,\ldots,e_m|e_ie_j+e_je_i=0,i\ne
j,e_i^2=1\rangle$. In $V_1$ we can find a basis
$v_1,w_1,\ldots,v_n,w_n$ such that $f(v_i,w_j)=\delta_{ij}$,
$f(v_i,v_j)=f(w_i,w_j)=0$. Consider the Weyl algebra $W_n=\langle
1,v_i,w_i,1\le i\le n,
[v_i,w_j]=\delta_{ij},[v_i,w_j]=[v_i,w_j]=0\rangle$.  According to
\cite{MZ}, the universal associative enveloping algebra of $F+V$
is isomorphic to the (super)tensor product $C(V,f)\otimes_F W_n$.
We will utilize this fact in the following Lemma.

\begin{lemma} There are no subsuperalgebras $\cal B$ isomorphic to
$J(V,f)$, where $V=V_0+V_1$, $V_1\ne\{0\}$ in a finite-dimensional
Jordan superalgebra ${\cal A}^{(+)}$, where $\cal A$ is an
associative superalgebra.
\end{lemma}
\begin{proof}
We assume the contrary, that is, there exists a subsuperalgebra
$\cal B$ of the type $J(V,f)$ in  ${\cal A}^{(+)}$. For $\cal B$,
we consider the universal associative enveloping superalgebra
$U({\cal B})$. According to the above fact, $U({\cal
B})=C(V_0,f)\otimes_F W_n$ where $C(V_0,f)$ is a Clifford algebra
for $V_0$, $f$ is a bilinear form on $V_0$, $W_n$ is a Weyl
algebra, $n=\frac{1}{2}\dim\,V_1$. Let $\varphi$ denote the
identity embedding of $\cal B$ in $\cal A$. As a direct
consequence of the definition of  universal enveloping algebra,
$\varphi$ can be uniquely extended to a homomorphism
$\bar\varphi:U({\cal B})\to {\cal A}$. Note that
$\bar\varphi(x)=\varphi(x)=x$ where $x\in V_1$. In other words,
$\bar\varphi(V_1)\ne 0$. However, since $V_1$ generates $W_n$,
$\bar\varphi(W_n)\ne 0$. It follows that $\bar\varphi(W_n)\cong
W_n$. Therefore, $\cal A$ has an infinite-dimensional
subsuperalgebra. This contradicts our assumptions.
\end{proof}

In the next Lemma we will prove that no simple decompositions in
which one of the components has either the type $K_3$ or $D_t$ are
possible.

\begin{lemma} Any superalgebra $\cal J$ of the type
$M_{n,m}(F)^{(+)}$, $n,m>0$ cannot be represented as the sum of
two proper simple subsuperalgebras  one of which has the type
$K_3$ or $D_t$.
\end{lemma}
\begin{proof} First of all, we note that if $\mbox{char}\,F=0$, then $K_3$
has no non-trivial finite-dimensional associative specializations
(see Remark 1). Therefore, we can directly pass to the second case
when one of the subsuperalgebras is isomorphic to $D_t$. Next we
suppose that $\cal J$ is given in the canonical form which is the
set of all matrices of order $(n+m)$ with respect to the natural
$Z_2$-gradation. Let $V=V_0+V_1$ denote the $Z_2$-graded module
corresponding to the natural representation of $M_{n,m}(F)^{(+)}$,
$\dim\,V_0=n$, $\dim\,V_1=m$. Any decomposition of
$M_{n,m}(F)^{(+)}$ induces that of $(M_{n,m}(F)^{(+)})_0=H({\cal
R}_n)\oplus H({\cal R}_m)$ given by
$$H({\cal R}_n)\oplus H({\cal R}_m)=(Fe_1\oplus Fe_2)+{\cal B}_0,$$
where $e_1$, $e_2$ are pairwise orthogonal idempotents in ${\cal
A}\cong D_t$. Next we define a pair of homomorphisms denoted as
$\pi_1$, $\pi_2$ which are the projections on the ideals $H({\cal
R}_n)$ and $H({\cal R}_m)$, respectively. Then the above
decomposition can be rewritten in the following way:
$$ H({\cal R}_n)=\pi_1(Fe_1\oplus Fe_2)+\pi_1({\cal B}_0),$$
$$ H({\cal R}_m)=\pi_2(Fe_1\oplus Fe_2)+\pi_2({\cal B}_0).$$
Next we estimate the dimension of $\pi_1({\cal B}_0)$. We know
that $\pi_1({\cal B}_0)$ is either a simple or a non-simple
semisimple subalgebra. Therefore, in the first case we have $\dim
\pi_1({\cal B}_0)\le n^2-2n+1$, and in the second case $\dim
\pi_1({\cal B}_0)\le n^2-2n+2$. It follows that $\dim H({\cal
R}_n)\le 2+\dim \pi_1({\cal B}_0)$, $n^2\le n^2-2n+4$, $n\le 2$.
Using the same arguments as above we can prove that $m\le 2$. As a
result, we have four possible cases ${\cal J}=M_{1,1}(F)^{(+)}$,
$M_{1,2}(F)^{(+)}$, $M_{2,1}(F)^{(+)}$ or $M_{2,2}(F)^{(+)}$.
Notice that the first case can be immediately excluded because
$\dim M_{1,1}(F)^{(+)}=4$. Since $M_{1,2}(F)^{(+)}\cong
M_{2,1}(F)^{(+)}$ the second and  forth are the only cases of
interest to us. Let ${\cal J}=M_{1,2}(F)^{(+)}$. By the dimension
argument, $\dim {\cal B}\ge 5$, and, moreover, $\mbox{rk}\,{\cal
B}\le 3$. Clearly, there are only three appropriate choices for
${\cal B}:$ $osp(2,1)$, $P(2)$ or $Q(2)$. However, for all cases,
$U({\cal B})$ is isomorphic to $M_{2,2}(F)$, that is, cannot be a
subsuperalgebra of $M_{1,2}(F)$.

Finally, let ${\cal J}=M_{2,2}(F)^{(+)}$. Again by the dimension
argument,  $\dim {\cal B}\ge 12$ and $\mbox{rk}\,{\cal B}_0\le 4.$
Considering all possible cases we come to the conclusion that
there are no appropriate subsuperalgebras in $\cal J$. This proves
our Lemma.

\end{proof}

\begin{lemma} Any superalgebra of the type $M_{n,m}(F)^{(+)}$, $n,m>0,$
cannot be represented as the sum of two proper non-trivial simple
subsuperalgebras one of which has either the types
$M_{1,1}(F)^{(+)}$, $osp(1,1)$ or $P(2)$.
\end{lemma}
\begin{proof}
Assume that $M_{n,m}(F)^{(+)}={\cal A}+{\cal B}$ where $\cal A$
and $\cal B$ satisfy all the above conditions. Let $\cal A$ have
one of types $M_{1,1}(F)^{(+)}$, $osp(1,1)$ or $P(2)$.  The even
part of the above decomposition can be rewritten as follows:
$$ H({\cal R}_n)=\pi_1({\cal A}_0)+\pi_1({\cal B}_0),$$
$$ H({\cal R}_m)=\pi_2({\cal A}_0)+\pi_2({\cal B}_0).$$

Let ${\cal A}$ be isomorphic to either $M_{1,1}(F)^{(+)}$ or
$osp(1,1)$. Then $\dim\pi_1({\cal A}_0)=2$. Acting in the same
manner as in Lemma 2.5, we obtain the following possibilities:
$n=m=1$, $n=1$, $m=2$ ($m=2$, $n=1$), $n=m=2$. Obviously, there
are no possible simple decompositions in the first case due to the
low dimension of $M_{1,1}(F)^{(+)}$. In the second and third cases
we have the following restrictions on the dimension  and the rank
of ${\cal B}_0$:  $\dim {\cal B}_0\ge 5$, $\mbox{rk}\,{\cal
B}_0\le 3$;  $\dim {\cal B}_0\ge 12$, $\mbox{rk}\, {\cal B}_0\le
4$. Considering all cases one after another we conclude that there
is no suitable choice for ${\cal B}_0$. Therefore,
$M_{n,m}(F)^{(+)}\ne {\cal A}+{\cal B}$.

 In the last case when ${\cal A}\cong P(2)$ there are the
following restrictions on indices: $n\le 3$, $m\le 3$. In other
words, $n=1$, $m=2$; $n=m=2$; $n=1$, $m=3$; $n=2$, $m=3$; $n=m=3$.
By the dimension and rank arguments there is no such ${\cal B}_0$.
The Lemma is proved.
\end{proof}

Next taking into account all previous Lemmas we list  simple
decompositions that might exist in $M_{n,m}(F)^{(+)}$. Let $\cal
A$ and $\cal B$ stand for the simple non-trivial Jordan
subsuperalgebras of $M_{n,m}(F)^{(+)}$.
\par\medskip
$$\qquad \left| \begin{array}{c|c|c} \hline
   & {\cal A}    & {\cal B}\\
  \hline
  1& M_{k,l}(F)^{(+)}    &  M_{p,q}(F)^{(+)}\\
  2& M_{k,l}(F)^{(+)}    &   P(q)\\
  3& M_{k,l}(F)^{(+)}    &   Q(p)\\
  4& P(k)            &   Q(l)\\
  5& P(k)            &   P(l)\\
  6& Q(k)            &   Q(l)\\
  7& osp(k,l)    &  M_{p,q}(F)^{(+)}\\
  8& osp(k,l)    &   Q(p)\\
  9& osp(k,l)    &   P(q)\\
 10& osp(k,l)    &  osp(p,q)\\
  \hline
  \end{array}\right|
$$
\par\medskip
Considering associative subalgebras $S({\cal A})$ and $S({\cal
B})$ generated by $\cal A$ and $\cal B$, respectively, we obtain a
new decomposition of the form $M_{n+m}(F)=S({\cal A})+S({\cal B})$
where $S({\cal A})$ and $S({\cal B})$ are associative subalgebras
of $M_{n+m}(F)$. Note that $S({\cal A})$ is a homomorphic image of
$U({\cal A})$. As a direct consequence of Theorem 2.2, $U({\cal
A})$ is either an associative simple algebra or a direct sum of
two or more simple pairwise isomorphic associative algebras.
\begin{lemma}
Let $\cal A$ be a proper non-trivial simple subsuperalgebra in
$M_{n,m}(F)^{(+)}$ where $n,m>0$. Then $S({\cal A})$ coincides
with $M_{n+m}(F)$ if and only if one of the following conditions
hold

(1) Either ${\cal A}\cong osp(p,q)$, $p+2q=n+m$, or

(2) ${\cal A}\cong P(n)$ for the case when $n=m$.
\end{lemma}
\begin{proof}
First, we note that the converse of this Lemma is obvious (see
Theorem 2.2). To prove that one of the above conditions holds in
the case when $S({\cal A})=M_{n,m}(F)$, then  we first show that
$\cal A$ cannot be of type $M_{k,l}(F)^{(+)}$ or $Q(p)$. If $\cal
A$ has the type $M_{k,l}(F)^{(+)}$, then $k+l<n+m$. By Theorem
2.2, $S({\cal A})$ is either a proper simple subalgebra of the
type $M_{k+l}(F)$ or a non-simple semisimple subalgebra of the
type $M_{k+l}(F)\oplus M_{k+l}(F)$. In both cases, $S({\cal A})\ne
M_{n+m}(F)$.

If ${\cal A}\cong Q(k)$, then its associative enveloping algebra
is a non-simple semisimple subalgebra which is the direct sum of
two or more simple ideals of the type $M_k(F)$. Therefore,
$S({\cal A})\ne M_{n+m}(F)$.

For the rest cases, $\cal A$ can either have  the type $osp(p,q)$
or $P(k)$. If ${\cal A}\cong osp(p,q)$, then $S({\cal A})\cong
M_{p+2q}(F)$. Hence $S({\cal A})=M_{n+m}(F)$ if and only if
$p+2q=n+m$. This yields (1).

Next we continue our proof by assuming that $n\ne m$, say, $n<m$.
We let $\cal A$ have the type $P(k)$. Then its even component
${\cal A}_0$, which is isomorphic to $H({\cal R}_k)$, is a proper
subalgebra in $ M_{n,m}(F)_{0}^{(+)}= I_1\oplus I_2$, $I_1\cong
H({\cal R}_n)$, $I_2\cong H({\cal R}_m)$. As previously, let
$\pi_1$ and $\pi_2$ denote the projections on $I_1$ and $I_2$,
respectively.

Suppose that $S({\cal A})=M_{n+m}(F)$. Since ${\cal A}\cong P(k)$,
then $S({\cal A})\cong M_{2k}(F)$. This implies $2k=n+m$,
$k=\frac{n+m}{2}$. In particular, $k>n$. Hence $\pi_1({\cal
A}_0)=\{0\}$ and $\pi_2({\cal A}_0)\cong H({\cal R}_k)$. It
follows that ${\cal A}_0\subseteq I_2$. Thus the identity $e$ of
$\cal A$ is an element of $I_2$. For any $x\in {\cal A}_1$,
$xe+ex=2x$ where the multiplication is associative. Multiplying
both sides of this equation by $e$, we obtain the following
$exe+ex=2ex$. Since $exe=0$, we have $ex=2ex$. Similarly, $xe=0$,
that is, $x=0$, for any $x\in {\cal A}_1$, a contradiction.

In conclusion, it remains to consider the case when $n=m$ and
${\cal A}\cong P(n)$. However,  it is obvious that $S({\cal
A})\cong M_{2k}(F)$ and $S({\cal A})=M_{2n}(F)$ if and only if
$k=n$. This completes our proof.
\end{proof}
\begin{lemma} Let $M_{n,m}(F)^{(+)}={\cal A}+{\cal B}$, $n,m>0$. Then one of the
subsuperalgebras in the given decomposition has either the type
$osp(p,q)$ where $p+2q=n+m$ or $P(n)$ (only if $n=m$).
\end{lemma}
\begin{proof}
Let us assume the contrary, that is, neither $\cal A$ nor $\cal B$
is a subsuperalgebra of any of the above types. Then, by Lemma
2.7, $S({\cal A})$ and $S({\cal B})$ are proper associative
subalgebras in $M_{n+m}(F)$. Theorem 2.2 states that both $S({\cal
A})$ and $S({\cal B})$ are either simple associative algebras or
non-simple semisimple associative algebras decomposable into the
sum of two or more pairwise isomorphic simple algebras. Therefore,
$\dim S({\cal A})\le k^2(\frac{n+m}{k})=(n+m)k$ where $k^2$ is a
dimension of a simple ideal,  $k>1$. If one of the
subsuperalgebras in the decomposition of $M_{n+m}(F)$ has a
non-zero annihilator then by Proposition in \cite{BK} no such
decomposition exists. Hence the identity of $M_{n+m}(F)$ is
contained in the intersection of $S({\cal A})$ and $S({\cal B})$.
On the other hand, the dimension of $S({\cal A})$ as well as
$S({\cal B})$ is strictly greater then $\frac{(n+m)^2}{2}$.

Thus, by the dimension argument, the sum of $S({\cal A})$ and
$S({\cal B})$ is a proper vector subspace of $M_{n+m}(F)$.
Therefore, $M_{n+m}(F)\ne S({\cal A})+S({\cal B})$. This implies
that our hypothesis was wrong.
\end{proof}

\begin{lemma}
Let $M_{n,m}(F)^{(+)}={\cal A}+{\cal B}$, $n,m>0$. Then, in the
case when $m$ is even, and $n$ is odd,    ${\cal A}\cong
osp(n,\frac{m}{2})$ and ${\cal B}\cong M_{k,l}(F)^{(+)}$ where
either $k=n-1,n$ or $l=m$. On the contrary, if $m$ is odd, and $n$
is even, then ${\cal A}\cong osp(m,\frac{n}{2})$ and ${\cal
B}\cong M_{k,l}(F)^{(+)}$ where either $k=m-1,m$ or $l=n$.
\end{lemma}
\begin{proof}

Since the proof remains the same for both cases, we consider only
the first case. First, let $n\ne m$. In view of Lemma 2.8, one of
the subsuperalgebras in $M_{n,m}(F)^{(+)}={\cal A}+{\cal B}$, for
example $\cal A$, is isomorphic to $osp(p,q)$  where
$$p+2q=n+m.\eqno (1) $$
The decomposition of $M_{n,m}(F)^{(+)}$ given above induces the
following representation of the even component
$M_{n,m}(F)_0^{(+)}={\cal A}_0+{\cal B}_0$ where
$M_{n,m}(F)_0^{(+)}=H({\cal R}_n)\oplus H({\cal R}_m)$, ${\cal
A}_0\cong H(F_p)\oplus H({\cal Q}_q)$. If for some $i$,
$\pi_i({\cal A}_0)\cong H(F_p)\oplus H({\cal Q}_q)$, then either
$p+2q\le n $ or $p+2q\le m$. However these inequalities conflict
with condition (1). Hence either $\pi_1({\cal A}_0)\cong H(F_p)$,
$\pi_2({\cal A}_0)\cong H({\cal Q}_q)$ or  $\pi_1({\cal A}_0)\cong
H({\cal Q}_q)$, $\pi_2({\cal A}_0)\cong H(F_p)$. If the first
possibility holds true, then

1. $H({\cal R}_n)=\pi_1({\cal A}_0)+\pi_1({\cal B}_0)$,
$\pi_1({\cal A}_0)\cong H(F_p)$, $p\le n$, $\pi_1({\cal B}_0)\ne
0$.

2. $H({\cal R}_m)=\pi_2({\cal A}_0)+\pi_2({\cal B}_0)$,
$\pi_2({\cal A}_0)\cong H({\cal Q}_q)$, $q\le\frac{m}{2}$,
$\pi_2({\cal B}_0)\ne 0$.

Since $p+2q=n+m$, it follows that $p=n$ and $q=\frac{m}{2}$.
Clearly, $\cal A$ has the type $osp(n,\frac{m}{2})$.
 If the second possibility holds true, then acting in the same manner, we
can show that $p=m$, $q=\frac{n}{2}$. However, we assumed that $n$
is odd. Hence it remains to prove that ${\cal B}\cong
M_{k,l}(F)^{(+)}$ where a pair of indices  $k,l$ satisfies the
conditions given in the Lemma. To prove this, we consider all
possible types for $\cal B$ in a step-by-step manner.

If ${\cal A}\cong osp(n,\frac{m}{2})$, ${\cal B}\cong P(k)$, then
the decomposition induces the following representation of the odd
part: $ M_{n,m}(F)^{(+)}_1={\cal A}_1+{\cal B}_1$ where $\dim
{\cal A}_1=nm$, $\dim {\cal B}_1=k^2$, that is, $2nm\le nm+k^2$,
$nm\le k^2$. Conversely,  $k\le n$, $k\le m$ since both
projections $\pi_1({\cal B}_0)$, $\pi_2({\cal B}_0)$ are non-zero.
Moreover, one of the inequalities should be strict since $n\ne m$.
Therefore, $k^2<nm$, which is a contradiction.

If  ${\cal A}\cong osp(n,\frac{m}{2})$, ${\cal B}\cong Q(k)$,
then, acting in the same manner as in the previous case, we can
prove that $M_{n,m}(F)^{(+)}\ne {\cal A}+{\cal B}$.

If ${\cal A}\cong osp(n,\frac{m}{2})$, ${\cal B}\cong osp(p,q)$,
then, by the dimension argument, we have the following inequality
$(n+m)^2\le 2(\frac{n(n+1)}{2}+\frac{m(m-1)}{2}+nm)$. Simplifying
the last inequality, we obtain that $m\le n$. Clearly, the
opposite inequality $m\ge n$ also holds true. Therefore, $m=n$
which contradicts our hypothesis.
 Overall, it remains to consider the case when ${\cal
A}\cong osp(n,\frac{m}{2})$, ${\cal B}\cong M_{k,l}(F)^{(+)}$.

Again the decomposition of $M_{n,m}(F)^{(+)}$ induces that of
$M_{n,m}(F)^{(+)}_0$ as follows: $M_{n,m}(F)^{(+)}_0={\cal
A}_0+{\cal B}_0$. Moreover, $M_{n,m}(F)^{(+)}_0= H({\cal
R}_n)\oplus H({\cal R}_m)$, ${\cal A}_0\cong H(F_n)\oplus H({\cal
Q}_{\frac{m}{2}})$, ${\cal B}_0\cong H({\cal R}_k)\oplus H({\cal
R}_l)$. If both $\pi_1({\cal B}_0)$ and $\pi_2({\cal B}_0)$ are
non-simple semisimple, that is, $\pi_1({\cal B}_0)\cong {\cal
B}_0$ and $\pi_2({\cal B}_0)\cong {\cal B}_0$, then we have the
following restrictions: $k+l\le n$ and $k+l\le m$. Since $n\ne m$,
we can assume without any loss of generality that $n<m$. Hence the
dimension of $\pi_i({\cal B}_0)$, $i=1,2$, is less than
$n^2-2n+2$. It follows from $\pi_1({\cal B}_0)\cong {\cal B}_0$
that $\dim\,{\cal B}_0\le n^2-2n+2$.

As a result, $\dim
M_{n,m}(F)^{(+)}_0=n^2+m^2\le\frac{n(n+1)}{2}+\frac{m(m-1)}{2}+n^2-2n+2$,
$\frac{m(m+1)}{2}\le\frac{n(n+1)}{2}+2-2n$, which is wrong.
Therefore, we have only two possibilities: either $\pi_1({\cal
B}_0)$ or $\pi_2({\cal B}_0)$ is a simple algebra. According to
\cite{tm}, for the first case, $k=n-1,n$ and, for the second,
$l=m$. Thus the Lemma is proved for the case when $n\ne m$.

To complete our proof we consider the case when $n=m$. First, we
assume that neither $\cal A$ nor $\cal B$ has the type $P(n)$. By
the previous Lemma one of the subsuperalgebras, for example $\cal
A$, is isomorphic to $osp(p,q)$, $p+2q=2n$, that is, $p=n$,
$q=\frac{n}{2}$. Then
$$ H({\cal R}_n)=\pi_1({\cal A}_0)+\pi_1({\cal B}_0),\quad \pi_1({\cal A}_0)\cong H(F_n), \pi_1({\cal B}_0)\ne 0$$
$$ H({\cal R}_n)=\pi_2({\cal A}_0)+\pi_2({\cal B}_0),\quad \pi_2({\cal A}_0)\cong H({\cal Q}_\frac{n}{2}),
\pi_2({\cal B}_0)\ne 0.$$

For some $i$, let $\pi_i({\cal B}_0)$ be a non-simple semisimple
subalgebra, then
$$\pi_i({\cal B}_0)\cong\left\{\begin{array}{c}
                H({\cal R}_k)\oplus H({\cal R}_l),\, k+l\le n
                \quad\mbox{or}\\
                 H(F_k)\oplus H({\cal Q}_l),\, k+2l\le n
                 \end{array}\right.$$

Therefore, $\dim\pi_i({\cal B}_0)\le n^2-2n+2.$ However $\dim
M_{n,n}(F)_0^{(+)}\le \dim{\cal A}_0+\dim{\cal B}_0$, $2n^2\le
n^2-2n+2+\frac{n(n+1)}{2}+\frac{n(n-1)}{2}=2n^2-2n+2$, that is,
$n\le 1$.  As mentioned above, there are no simple decompositions
in $M_{1,1}(F)^{(+)}$. Hence both $\pi_1({\cal B}_0)$ and
$\pi_2({\cal B}_0)$ are simple. It follows that $\pi_1({\cal
B}_0)\cong H({\cal R}_{n-1})$, $\pi_2({\cal B}_0)\cong H({\cal
R}_n)$, that is, ${\cal B}_0\cong M_{n-1,n}(F)$.

Next we let $\cal A$ be of the type $P(n)$. Then ${\cal
B}\cong\left\{\begin{array}{c}
                                 P(k)\\
                                 Q(k)\\
                                 osp(k,l)\\
                                 M_{k,l}(F)
                                \end{array}\right.,$
for some integers $k$ and $l$.

{\bf 1.} ${\cal B}\cong P(k)$, hence $\dim\,{\cal B}=2k^2$, $k\le
n$. For $\dim M_{n,n}(F)^{(+)}\le \dim {\cal A}+\dim {\cal B}$, it
is clear that $k=n$, and the sum in the decomposition is direct.
However since both subsuperalgebras have the type $P(n)$, they
contain the identity of $M_{n,n}(F)^{(+)}$, a contradiction.

{\bf 2.} ${\cal B}\cong Q(k)$. In this case the proof is the same
as in Case 1.

{\bf 3.} ${\cal B}\cong osp(k,l)$. Clearly, for some $i$,
$\pi_i({\cal B}_0)$ is non-simple semisimple. Therefore, $k+2l\le
n$. In particular,  $\pi_i({\cal B}_0)\cong {\cal B}_0$. Hence,
$\dim{\cal B}_0\le n^2-2n+2$. Thus $\dim
M_{n,n}(F)_0^{(+)}=2n^2\le 2n^2-2n+2$, a contradiction. As a
result, $\pi_i({\cal B}_0)$, $i=1,2$, is a simple subalgebra, that
is, $k\le n$, $l\le\frac{n}{2}$. Then $\dim {\cal B}\le 2n^2$. By
the dimension argument, $k=n$, $l=\frac{n}{2}$ and the sum in the
given decomposition is direct. However, this contradicts the fact
that both subsuperalgebras in the given decomposition contain the
identity of $M_{n,n}(F)^{(+)}$.

{\bf 4.} ${\cal B}\cong M_{k,l}(F)^{(+)}$, $k+l<2n$. The even part
of $M_{n,n}(F)^{(+)}$, that is, $M_{n,n}(F)_0^{(+)}$ equals to the
sum of two orthogonal ideals denoted as $I_1$ and $I_2$, both
ideals isomorphic to $H({\cal R}_n)$. By the dimension argument,
$\dim M_{n,n}(F)^{(+)}\le 2n^2+(k+l)^2$, $4n^2\le 2n^2+(k+l)^2$,
$k+l\ge\sqrt 2 n$. In particular, $\pi_1({\cal B}_0)$,
$\pi_2({\cal B}_0)$ are simple. Therefore, acting by a appropriate
automorphism of $M_{n,n}(F)^{(+)}$, ${\cal B}_0$ can be reduced to
the block-diagonal form. Moreover, $I_1$ and $I_2$ contain all
simple ideals isomorphic to $H({\cal R}_k)$ and $H({\cal R}_l)$,
respectively.

Suppose that the identity of $M_{n,n}(F)^{(+)}$ is an element of
$\cal B$. This implies that $kk_1=ll_1=n$ where $k_1$ and $l_1$
are the numbers of blocks which have types $H({\cal R}_k)$ and
$H({\cal R}_l)$, respectively. In view of the inequality
$k+l\ge\sqrt 2 n$ this result implies that either $k_1=2$, $l_1=1$
or $k_1=1$, $l_1=2$, that is, ${\cal B}\cong M_{\frac{n}{2},n}(F)$
up to the order of indices. By Theorem 2.2, $S({\cal B})$ is
isomorphic to $M_{\frac{3n}{2}}(F)$ or $M_{\frac{3n}{2}}(F)\oplus
M_{\frac{3n}{2}}(F)$. Obviously, $S({\cal B})$ cannot be
non-simple semisimple of the indicated type because its rank is
greater than $2n$. However, by Lemma 1.1, the first case is also
impossible because $\frac{3n}{2}>n$. Hence the identity of
$M_{n,n}(F)$ is not an element of $\cal B$. In other words, $\cal
B$ as well as ${\cal B}_0$ has a non-zero annihilator.

 Acting by appropriate automorphism of $M_{n,n}(F)$ we can reduce $\cal
B$ to the following form:
 $$\left\{ \left(\begin{array}{c|ccc|ccc}
0      &   0   &  \ldots  &   0   &  0  & \ldots   &  0\\
\hline
0      &       &          &       &     &          &    \\
\vdots &       &    T_1   &       &     &   T_2    &    \\
0      &       &          &       &     &          &    \\
\hline
0      &       &          &       &     &          &    \\
\vdots &       &    T_3   &       &     &   T_4    &    \\
0      &       &          &       &     &          &
\end{array}\right)\right\}
$$
where $T_1$, $T_2$, $T_3$ and $T_4$ are matrices of orders
$(n-1)\times (n-1)$, $(n-1)\times m$, $m\times(n-1)$ and $m\times
m$, respectively.

This implies that ${\cal A}_0$ takes the form:
$$
\left\{\left(\begin{array}{cc}
            X& 0\\
            0& C^{-1}X^tC
            \end{array}\right)\right\},
$$
for some $C$, $\mbox{det}\,C\ne 0$. Then using the automorphism
$\varphi(Y)=C'^{-1}YC'$ where
$$     C'=\left(\begin{array}{cc}
            I& 0\\
            0& C^{-1}
            \end{array}\right),
$$
$\cal A$ can be reduced to the form where
$${\cal A}_0=\left\{\left(\begin{array}{cc}
                     X& 0\\
                     0& X^t
                    \end{array}\right)\right\}
$$
while $\cal B$ remains the same. Obviously, this decomposition is
not possible. The Lemma is proved.
\end{proof}

\begin{fact}
A Jordan superalgebra of the type $M_{n,m}(F)^{(+)}$ where $m$ is
even can be represented as the sum of two proper simple
subsuperalgebras $\cal A$ and $\cal B$ which have types
$osp(n,\frac{m}{2})$ and $M_{n-1,m}(F)^{(+)}$, respectively.
\end{fact}
\begin{proof}
To prove, we consider the first subsuperalgebra in the standard
realization:
$$\left\{ \left(\begin{array}{c|c}
                  A   & C\\
                \hline
           S^{-1}C^t  & B
           \end{array}\right)\right\}
$$
where $A$ is a symmetric matrix of order $n$, $B$ is a symplectic
matrix of order $m$, $C$ is any matrix of order $n\times m$,
$S=\left(\begin{array}{cc}
                  0& I\\
                  -I&0
               \end{array}\right)$
where $I$ is identity matrix of order $\frac{m}{2}$. The second
subalgebra can be viewed in the following form:
$$\left\{ \left(\begin{array}{ccc|cccc}
        &        &   & &    & &\\
        & A      &   & &  B & & \\
        &        &   & &    & &  \\
       \hline
       0&\ldots  & 0 & 0&\ldots   &0&\\
       \hline
        &        &   & &    & &\\
        &    C   &   & &   D &  & \\
        &        &   & &    &  &
        \end{array}\right)\right\}
$$
where $A$ and $C$ of orders $(n-1)\times n$ and $m\times n$,
respectively, have the last two columns equal, $B$ and $D$ are any
matrices of orders $(n-1)\times n$, $m\times m$, respectively. By
straightforward calculations $\dim ({\cal A}_1+{\cal B}_1)=\dim
{\cal A}_1+\dim {\cal B}_1-\dim({\cal A}_1\cap {\cal
B}_1)=mn+2m(n-1)-m(n-2)=2mn$. This proves our Lemma.
\end{proof}

\begin{fact}
In the case when $n$ is even, a Jordan algebra of the type
$M_{n,m}(F)^{(+)}$ can also be decomposed into the sum of $\cal A$
and $\cal B$ where ${\cal A}\cong osp(m,\frac{n}{2})$ and ${\cal
B}\cong M_{m-1,n}(F)^{(+)}$. This decomposition can be constructed
in the same manner as in Example 1.
\end{fact}

\begin{Prop}
Example 1 and 2 are  the only possible decompositions of
$M_{n,m}(F)^{(+)}$, $n,m>0$ into the sum of two proper simple
non-trivial subsuperalgebras for appropriate values of $n$, $m$.
\end{Prop}
\begin{proof}
As usual, we assume the contrary, that is, there exists some other
simple decomposition of $M_{n,m}(F)^{(+)}$ different from one in
Example 1. By Lemma 2.9, this decomposition takes the following
form:

1. If $m$ is even, then  $M_{n,m}(F)^{(+)}={\cal A}+{\cal B}$,
${\cal A}\cong osp(n,\frac{m}{2})$, ${\cal B}\cong
M_{l,k}(F)^{(+)}$ where either  $l=n-1, n$ or $k=m$.

2. If $n$ is even,  then  $M_{n,m}(F)^{(+)}={\cal A}+{\cal B}$,
${\cal A}\cong osp(m,\frac{n}{2})$, ${\cal B}\cong
M_{k,l}(F)^{(+)}$ where either $k=m-1,m$ or $l=n$.

Then $M_{n,m}(F)_1={\cal A}_1+{\cal B}_1$. It follows that $\dim
M_{n,m}(F)_1\le\dim{\cal A}_1+\dim{\cal B}_1$, that is, $2nm\le
nm+2lk$, $nm\le 2lk$. Hence, for even $m$,  $l\ge\frac{n}{2}$, in
the case $k=m$, and  $k\ge\frac{m}{2}$,  in the case $l=n-1$ or
$n$. Likewise, if $n$ is even, then $k\ge\frac{m}{2}$, in the case
$l=n$, and  $l\ge\frac{n}{2}$,  in the case $k=m-1$ or $m$.   For
definitness, we consider the case when $m$ is even, and  $l=n-1$
because the proof remains the same for all other cases.

Let $V=V_0+V_1$ denote a $Z_2$-graded vector space where
$\dim\,V_0=n$ and $\dim\,V_1=m$. By its definition,
$M_{n,m}(F)^{(+)}$ coincides with the set of all linear
transformations acting in $V$. Then let $\rho$ stand for the
natural representation of ${\cal B}={\cal B}_0+{\cal B}_1$  in
$V$. It follows from the definition of this action that
$\rho({\cal B}_0)(V_0)\subseteq V_0$, $\rho({\cal
B}_0)(V_1)\subseteq V_1$, $\rho({\cal B}_1)(V_0)\subseteq V_1$,
$\rho({\cal B}_1)(V_1)\subseteq V_0$. Since $\cal B$ is a
non-simple semisimple Jordan algebra it acts completely reducibly
in $V$. Next we describe this action in more details. For this, we
identify $V$ with a $Z_2$-graded  vector space of the form
$W=\langle v_0\rangle \oplus (V_0'\otimes F^{r+1})\oplus V_1'$,
$r\ge 1$ where $v_0$ is a vector in $V_0$ annihilated by ${\cal
B}_0$, $V_0'$ is an invariant complementary subspace of $\langle
v_0\rangle$, $\rho({\cal B}_0)|_{V_0'}\cong H({\cal R}_{n-1})$,
$V_1'$ is an invariant subspace of $V_1$ such that  ${\cal B}_0$,
$\rho({\cal B}_0)|_{V_1'}\cong H({\cal R}_k)$. Moreover,
$W_0=\langle v_0\rangle \oplus V_0'\otimes e_0$, $W_1=V_0'\otimes
\langle e_1,\ldots,e_r\rangle\oplus V_1'$ where $\langle
e_0,e_1,\ldots,e_r\rangle$ is a basis for $F^{r+1}$. Then,
$\rho({\cal B}_0)=\rho({\cal B}_0)|_{\langle v_0\rangle}\oplus
\rho({\cal B}_0)|_{V_0'}\otimes Id_{r+1}\oplus\rho({\cal
B}_0)|_{V_1'}.$ Note that $\rho({\cal B}_0)|_{\langle
v_0\rangle}=0$. In other words, by choosing an appropriate  basis
in $V_0$ and $V_1$, $\rho({\cal B}_0)$ can be written in a
block-diagonal form in which the first block of order 1 is zero,
the last block has order $k$, and the other blocks have order
$r+1$. Next we consider the representation of the odd part ${\cal
B}_1$. For this, we choose any $a\in {\cal B}_0$ such that
$$\rho(a)(V_0'\otimes F^{r+1})=0,
\quad \rho(a)(V_1')\ne 0.\eqno (2)$$

All such elements form an ideal of ${\cal B}_0$ isomorphic to
$H({\cal R}_k)$. Then we choose any non-zero $x$ in ${\cal B}_1$.
Let $e$ denote the identity of $\cal B$, $e\in {\cal B}_0$. Then
$\rho(x)v_0=\rho(x\odot
e)v_0=\rho(\frac{xe+ex}{2})v_0=\frac{1}{2}(\rho(x)\rho(e)v_0+\rho(e)\rho(x)v_0)=
\frac{1}{2}\rho(x)v_0$, that is, $\rho(x)v_0=0$, for any $x\in
{\cal B}_1$. Next we find the representation of $a\odot x\in{\cal
B}_1$. As mentioned above, $\rho(a\odot x)(v_0)=0$. Besides,
$2\rho(a\odot x)(V_0'\otimes e_0)=\rho(a)\rho(x)(V_0'\otimes
e_0)+\rho(x)\rho(a)(V_0'\otimes e_0)\subseteq V_1'$, $\rho(a\odot
x)(V_0'\otimes \langle e_1,\ldots,e_r\rangle)=0$, $\rho(a\odot
x)(V_1')\subseteq V_0'\otimes e_0$.  Clearly, we can  find
$c\in{\cal B}_0$ whose action is given by the following formulae:
$$\rho(c)(V_0'\otimes F^{r+1} )\ne 0,\quad
\rho(c)(V_1')=0. \eqno (3)$$

Now we need to determine
$$ c\odot(x\odot a). \eqno (4)$$
Since $2\rho(c\odot (a\odot x))=\rho(c)\rho(a\odot x)+\rho(a\odot
x)\rho(c),$ we have the following: $\rho(c\odot (a\odot
x))(v_0)=0$, $\rho(c\odot (a\odot x))(V_0'\otimes \langle
e_1,\ldots, e_r\rangle)=0$. Besides,
$$\rho(c\odot (x\odot a))(V_0'\otimes e_0)=\rho(c)\rho(x)\rho(a)(V_0'\otimes e_0)+\rho(x)\rho(c)\rho(a)(V_0'\otimes
e_0)+$$
 $$\rho(a)\rho(c)\rho(x)(V_0'\otimes e_0)+\rho(a)\rho(x)\rho(c)(V_0'\otimes e_0)=
\rho(a)\rho(x)\rho(c)(V_0'\otimes e_0)\subseteq V_1'. \eqno (5)$$
Similarly,
$$\rho(c\odot (x\odot a))(V_1')=\rho(c)\rho(x)\rho(a)(V_1')+\rho(x)\rho(c)\rho(a)(V_1')
+\rho(a)\rho(c)\rho(x)(V_1')+$$ $$\rho(a)\rho(x)\rho(c)(V_1')=
\rho(c)\rho(x)\rho(a)(V_1')\subseteq V_0'\otimes e_0.\eqno (6)$$
Assume that $\rho(x)(V_1')\ne 0$, $\rho(x)(V_0'\otimes e_0)\ne
0(\mbox{mod}\,V_0'\otimes\langle e_1,\ldots,e_r\rangle)$. Then
$\rho(c\odot (x\odot a))$ has the following matrix form:
$$\left( \begin{array}{c|c|cc}
          0&\ldots & 0\ldots & 0\\
    \hline
   \vdots  & 0    & 0\ldots        & XY_1Z\\
    \hline
    0  &    0   &         & \\
\vdots &  \vdots     &         &\\
      0& ZY_2X &         & 0
    \end{array}\right),
 $$
where $X$ is an arbitrary square matrix of order $k$, $Y_1$ and
$Y_2$ are some fixed non-zero matrices of order $k\times (n-1)$
and $(n-1)\times k$, respectively, $Z$ is any square matrix of
order $n-1$. Next we choose any $y\in {\cal B}_1$. We have seen
that there exists an element $y$ of form (4) such that
$\rho(y-a\odot (x\odot c))(V_1')=0$ or $\rho(y-a\odot (x\odot
c))(V_0'\otimes e_0)=0(\mbox{mod}\, V_0'\otimes \langle
e_1,\ldots,e_r\rangle)$. Suppose that one of the above equations
does not hold. Without any loss of generality we let
$\rho(y')(V_0'\otimes e_0)\ne 0(\mbox{mod}\, V_0'\otimes\langle
e_1,\ldots,e_r\rangle)$, where $y'=y-a\odot (x\odot c)$.
Multiplying $y'$ by the elements of the form (2) and then (3) we
obtain $a'\odot(y'\odot c')\in{\cal B}_1$, where $a'$ and $c'$ ran
 relevant sets, and $\rho(a'\odot (y'\odot c'))(V_0'\otimes
F^{r+1})=0$, $\rho(a'\odot (y'\odot
c'))V_1'=\rho(a')\rho(y')\rho(c')V_1'\subseteq V_0'\otimes e_0$.
Moreover, $\rho(a'\odot(y'\odot c')):V_1'\to V_0'\otimes e_0$
represents all linear transformations from $k$-dimensional vector
space into $(n-1)$-dimensional vector space. Besides, all such
elements are linearly independent from all the elements (4).
Therefore, we found $2(n-1)k$ linearly independent elements of
${\cal B}_1$, ($\dim {\cal B}_1=2(n-1)k$). If there is at least
one element $\bar y\in {\cal B}_1$ such that $\rho(\bar
y)(V_0'\otimes e_0)\ne 0(\mbox{mod}V_1')$ or $\rho(\bar
y)(V_0'\otimes\langle e_1,\ldots,e_r\rangle)\ne 0$, then it will
be also linearly independent with all above elements. Hence, by
dimension arguments, there is no $\bar y$ satisfying the above
conditions. Consequently, for all elements in ${\cal A}_1$, the
following $\rho(\bar y)(V_0')=0(\mbox{mod} V_1')$, $\rho(\bar
y)(V_0'\otimes\langle e_1,\ldots,e_r\rangle)=0$, $\rho(\bar
y)(V_1')\subseteq V_0'\otimes e_0$ hold true. If we fix a basis in
$V$ such that in this basis the even part has the diagonal form:
$$ \left(\begin{array}{c|ccc}
            X&    &   0    &  \\
            \hline
             &   X& \ldots & 0 \\
            0&   \vdots & \ddots &   \\
             &   0 &        & Y
           \end{array}\right),
$$
then the odd part becomes the following:
$$  \left(\begin{array}{c|ccc}
                   0&  0& \ldots 0& Z\\
           \hline
                   0&   &        &  \\
           \vdots   &   &    0   &  \\
                   0&   &        & \\
                  Z'&   &        &
           \end{array}\right)\eqno (7)
$$
\noindent where $Z$, $Z'$ are any matrices of order $(n-1)\times
k$ and $k\times (n-1)$, respectively. Then it follows from ${\cal
B}_1\odot {\cal B}_1\subseteq {\cal B}_0$ that ${\cal B}_1=0$, a
contradiction.

We henceforth assume that the equations $\rho(y-a\odot (x\odot
c))(V_1')=0$ and $\rho(y-a\odot (x\odot c))(V_0'\otimes
e_0)=0\,(\mbox{mod}\, V_0'\otimes\langle e_1,\ldots,e_r\rangle)$
hold true simultaneously. Then multiplying $y-a\odot(x\odot c)$ by
the elements (4), we obtain  some elements of ${\cal B}_0$ which
act on $V_1'$ and $V_0'\otimes\langle e_1,\ldots,e_r\rangle$
non-invariantly. Hence, $y-a\odot(x\odot c)=0$. Therefore, the odd
component of ${\cal A}_1$ has form (7). As proved before, this is
not possible.

Next we assume that $\rho(x)(V_0'\otimes
e_0)=0(\mbox{mod}V_0'\otimes\langle e_1,\ldots,e_r\rangle)$, for
all $x\in {\cal B}_1$, and for at least one element $x'\in {\cal
B}_1$, $\rho(x')(V_1')\ne 0$.

Acting in the same manner as before, we obtain $a'\odot (x'\odot
c')\in {\cal B}_1$ which acts trivially on all subspaces except
for $V_1'$, which it carries into $V_0'\otimes e_0$. Considering
the difference between an arbitrary element $y\in {\cal B}_1$ and
a corresponding element $a''\odot (x''\odot c'')$, we can show
that $\rho(y-a''\odot(x''\odot c''))(V_0'\otimes
e_0)=0(\mbox{mod}\,V_0'\otimes\langle e_1,\ldots,e_r\rangle\oplus
V_1')$, $\rho(y-a''\odot(x''\odot c''))(V_1')=0$. Again
multiplying $a'\odot (x'\odot c')$ and $y-a''\odot (x''\odot
c'')$, we obtain some elements from ${\cal B}_0$ acting on $V_1'$
non-trivially. Then we conclude that ${\cal B}_1$ consists of all
elements which act on  $V_0'\otimes e_0$ trivially and carry the
other subspaces into $V_0'\otimes e_0$. Hence ${\cal B}_1\odot
{\cal B}_1=0$, a contradiction.

Finally, if $\rho(x)(V_1')=0$, $\rho(x)(V_0'\otimes
e_0)=0(\mbox{mod}\,(V_0'\otimes F^{r+1}))$, then it follows that
${\cal B}_1\odot {\cal B}_0=0$, which is clearly a wrong
statement. The Proposition is proved.
\end{proof}

Based on all above Lemmas and Proposition 1, we conclude that
Theorem 1 is true. In other words, $M_{n,m}(F)^{(+)}$ where
$n,m>0$, $m$ is even, and $n$ is odd admits only one decomposition
into the sum of two proper simple subsuperalgebras. If both
$n,\,m$ are odd, then $M_{n,m}(F)^{(+)}$ cannot be represented as
the sum of two proper simple subsuperalgebras. If both indices are
even, then $\cal A$ admits two different types of decompositions
of the following forms:

1. ${\cal A}={\cal B}_1+{\cal C}_1$ where ${\cal B}_1$ and ${\cal
C}_1$ have types $osp(n,\frac{m}{2})$ and $M_{n-1,m}(F)^{(+)}$,

2. ${\cal A}={\cal B}_2+{\cal C}_2$ where ${\cal B}_2$ and ${\cal
C}_2$ have types $osp(m,\frac{n}{2})$ and $M_{m-1,n}(F)^{(+)}$.

\section{Decompositions of  superalgebras of the type $osp(n,m)$}

This section is dedicated to the study of  simple decompositions
of $osp(n,m)$. Actually, we will show that  there are no such
decompositions over algebraically closed field $F$ of zero
characteristic. Our main purpose is to prove the following.

\begin{theorem}
Let $\cal J$ be a superalgebra of the type $osp(n,m)$ where
$n,m>0$. Then $\cal J$ cannot be written as the sum of two proper
nontrivial simple subsuperalgebras $\cal A$ and $\cal B$.
\end{theorem}

The proof of this Theorem is based on the following Lemmas.
\begin{lemma}
Let $\cal J$ be a superalgebra of type $osp(n,m)$ where $n,m>0$,
and $\cal A$, $\cal B$ are two proper simple subsuperalgebras none
of which has any of the types $K_3$ or $D_t$. Then $\cal J$ cannot
be represented as the sum of $\cal A$ and $\cal B$.
\end{lemma}

\begin{proof}
First  we identify $\cal J$ with $osp(n,m)$ which can be
considered in the canonical form. Next we assume the contrary,
that is,
$$osp(n,m)= {\cal A}+ {\cal B},\eqno (8)$$
The decomposition (8) generates the following decomposition of the
associative enveloping algebra into the sum of three non-zero
subspaces.
$$M_{n+2m}(F)= S(osp(n,m))=S({\cal A}) + S({\cal B})+S({\cal A})S({\cal
B}),\eqno (9)
$$
\noindent where $S({\cal A})$, $S({\cal B})$ denote the
associative enveloping algebras of $\cal A$, $\cal B$,
respectively. Let 1 denote the identity of $osp(n,m)$. Then we
consider the following cases.

{\bf Case 1.} Let $1\notin {\cal A}$, $1\notin {\cal B}$. This
implies that there exist non-zero $a_0$ and $b_0$ in
$\mbox{Ann}({\cal A})$ and $\mbox{Ann}({\cal B})$, respectively.
Then multiplying every term of (9) by $a_0$ on the left and $b_0$
on the right, the following equation $a_0M_{n+2m}(F)b_0=0$ takes
place, which is clearly wrong.

{\bf Case 2.} $1\in {\cal A}$, $1\in {\cal B}$. Six cases arise:

{\bf(a)} ${\cal A}\cong M_{k,l}(F)^{(+)}$, ${\cal B}\cong
M_{p,q}(F)^{(+)}$. The given decomposition induces the
decomposition of the even part $osp(n,m)_0={\cal A}_0+{\cal B}_0$
which in turn can be projected on the ideals of the even
component. In particular, $H(F_n)=\pi_1({\cal A}_0)+\pi_1({\cal
B}_0)$. By Theorem 1.3,  both projections cannot be simultaneously
simple. Therefore, at least one of the components is non-simple
semisimple. For definiteness, let $\pi_1({\cal A}_0) \cong H({\cal
R}_k)\oplus H({\cal R}_l)$. By Lemma  1.1, $k+l\le\frac{n}{2}$.
Then,  $\dim\pi_1({\cal A}_0)=k^2+l^2\le
(\frac{n}{2}-1)^2+1=\frac{n^2}{4}-n+2$. If $\pi_1({\cal B}_0)\cong
H({\cal R}_p)$  or $H({\cal R}_q)$, then, by Lemma 3.3,
$p\le\frac{n}{2}$ or $q\le\frac{n}{2}$. If $\pi_1({\cal B}_0)\cong
H({\cal R}_p)\oplus H({\cal R}_q)$, then $\dim\pi_1({\cal
B}_0)=p^2+q^2\le(\frac{n}{2}-1)^2+1=\frac{n^2}{4}-n+2$. As a
result, $\dim H(F_n)=\frac{n^2+n}{2}\le 2(\frac{n^2}{4}-n+2)$,
$\frac{5n}{2}\le 4$, $n\le 1$ or  $\dim H(F_n)=\frac{n^2+n}{2}\le
\frac{n^2}{4}-n+2+\frac{n^2}{4}$, $\frac{3n}{2}\le 2$, $n\le 1$.

There remains one case where $n=1$. The decomposition takes the
following form: $osp(1,m)={\cal A}+{\cal B}$ where ${\cal A}\cong
M_{1,l}(F)^{(+)}$, ${\cal B}\cong M_{p,q}(F)^{(+)}$. Then $H({\cal
Q}_m)=\pi_2({\cal A}_0)+\pi_2({\cal B}_0)$. If $\pi_2({\cal
A}_0)\cong F\oplus H({\cal R}_l)$ и $\pi_2({\cal B}_0)\cong
H({\cal R}_p)\oplus H({\cal R}_q)$, then $1+l\le m$, $p+q\le m$,
$\dim\pi_2({\cal A}_0)\le m^2-2m+2$, $\dim\pi_2({\cal B}_0)\le
m^2-2m+2$. As a result, $2m^2-4m+4\ge m(2m-1)=2m^2-m$, $4\ge 3m$,
$m\le 1$, that is, $m=1$. However, it is clear that both
subsuperalgebras in $osp(1,1)={\cal A}+{\cal B}$ are isomorphic to
$M(1,1)^{(+)}$, and, by dimension argument, $\cal A$, $\cal B$
coincide with $M(1,1)^{(+)}$.

If one of the projections is non-simple semisimple and the other
is simple, then $m(2m-1)\le m^2+m^2-2m+2=2m^2-2m+2$, $m\le 2$.
Therefore, $osp(1,2)={\cal A}+{\cal B}$, ${\cal A}\cong
M_{1,1}(F)^{(+)}$, ${\cal B}\cong M_{1,2}(F)^{(+)}$, $H({\cal
Q}_2)=\pi_2({\cal A}_0)+\pi_2({\cal B}_0)$ where $\pi_2({\cal
A}_0)$ is simple, and $\pi_2({\cal B}_0)$ is non-simple. By the
dimension argument, the sum in the above decomposition is direct,
that is, one of the subalgebras does not contain the identity,
which is obviously wrong. If both projections are simple, then
$osp(1,m)={\cal A}+{\cal B}$, $\cal A$, ${\cal B}\cong
M(1,m)^{(+)}$. However, by the dimension argument, the latter does
not hold.  Otherwise, ${\cal A}_1={\cal B}_1=osp(1,m)_1$ because
their dimensions are equal. It follows that ${\cal
A}_0=\langle{\cal A}_1\odot {\cal A}_1\rangle=\langle{\cal
B}_1\odot {\cal B}_1\rangle={\cal B}_0$, that is, ${\cal A}={\cal
B}=osp(1,m)$.

{\bf (2)}${\cal A}\cong M_{k,l}(F)^{(+)}$, ${\cal B}\cong P(q)$ or
$Q(q)$ ($q>1$). Therefore, $H(F_n)=\pi_1({\cal A}_0)+\pi_1({\cal
B}_0)$ where $\pi_1({\cal A}_0)\cong H({\cal R}_k)\oplus H({\cal
R}_l)$, $\pi_1({\cal B}_0)\cong H({\cal R}_q)$. Again, by the same
arguments as in the previous case, $n\le 1$. Let $n=1$ then
$\pi_1({\cal A}_0)\cong F$. In this case, $k=1$ (or $l=1$) and the
following decomposition holds true: $ H({\cal Q}_m)=\pi_2({\cal
A}_0)+\pi_2({\cal B}_0)$ where $\pi_2({\cal A}_0)\cong F\oplus
H({\cal R}_l)$, $\pi_2({\cal B}_0)\cong H({\cal R}_q)$ or
$\pi_2({\cal A}_0)\cong H({\cal R}_l)$, $\pi_2({\cal B}_0)\cong
H({\cal R}_q)$.

In the first case, we have proved that $m=2$, $l=1$. Hence
$osp(1,2)={\cal A}+{\cal B}$ where ${\cal A}\cong
M_{1,1}(F)^{(+)}$, ${\cal B}\cong P(2)$, which induces the
following:  $H({\cal Q}_2)=\pi_2({\cal A}_0)+\pi_2({\cal B}_0)$
where ${\cal A}_0\cong F\oplus F$, ${\cal B}_0\cong H({\cal
R}_2)$. The sum in the last decomposition is direct, and both
subalgebras contain 1, which is a contradiction. In the second
case, $1=q=m$, that is, $osp(1,m)={\cal A}+{\cal B}$, ${\cal
A}\cong M(1,m)^{(+)}$, ${\cal B}\cong P(m)$, $m>1$.

Since $\dim\,osp(1,m)_1\ge \dim {\cal B}_1$, then $2m\ge m^2$,
that is,  $m=2$. If $m=2$, then $osp(1,2)={\cal A}+{\cal B}$ where
${\cal A}\cong M_{1,2}(F)^{(+)}$, ${\cal B}\cong P(2)$ which
induces the equality  $H({\cal Q}_2)=\pi_2({\cal A}_0)+\pi_2({\cal
B}_0)$ where ${\cal A}_0\cong H({\cal R}_2)$, ${\cal B}_0\cong
H({\cal R}_2)$.

Notice that the identity of $osp(1,2)$ is an element of  $\cal A$,
that is, $\cal A$ has trivial two-sided annihilator. Consider an
associative enveloping algebras of  $osp(1,2)$ and $\cal A$
denoted as $S(osp(1,2))$ and $S({\cal A})$, respectively. It can
be shown that $S({\cal A})\cong M_3(F)$ is a subalgebra of
$S(osp(1,2))=M_5(F)$ (see Theorem 2.2). By Lemma 1.1, $S({\cal
A})$ contains no identity of $M_5(F)$, therefore, has a non-zero
two-sided annihilator, and so does $\cal A$, a contradiction.

{\bf(c)} ${\cal A}$, $\cal B$ have types $P(q)$ or $Q(p)$. Then
the decomposition leads to the decomposition of $H(F_n)$ into the
sum of two proper subalgebras, which does not exist.

{\bf (d)}  ${\cal A}\cong osp(k,l)$, ${\cal B}\cong
M_{p,q}(F)^{(+)}$. Since $S({\cal A}) \cong M_{k+2l}(F)$ contains
the identity of the entire superalgebra, $k+2l\le\frac{n+2m}{2}$.
Similarly,  $p+q\le \frac{n+2m}{2}$.

Thus $\dim osp(n,m)\le \dim {\cal A}+\dim {\cal B}$, that is,
$\frac{n^2+n}{2}+m(2m-1)+2nm\le
\frac{k^2+k}{2}+l(2l-1)+2kl+\frac{(n+2m)^2}{4}$. By
straightforward calculations we obtain
$\frac{n^2}{4}+\frac{n}{2}+m^2+nm\le 3m$, which is true if and
only if $m=n=1$. Obviously, $osp(1,1)$ has no simple
decompositions.

{\bf (e)} ${\cal A}\cong osp(k,l)$, ${\cal B}\cong P(q)$. Then, we
have  $k+2l\le\frac{n+2m}{2}$, $2q\le\frac{n+2m}{2}$. Therefore,
$\dim {\cal B}=2q^2\le 2(\frac{n+2m}{4})^2$. Again, by the
dimension argument, this decomposition is not possible.

{\bf (f)} ${\cal A}\cong osp(k,l)$, ${\cal B}\cong osp(p,q)$. Then
$k+2l\le\frac{n+2m}{2}$, $p+2q\le\frac{n+2m}{2}$. Comparing $\dim
osp(n,m)$ with $\dim\, {\cal A}+ \dim\, {\cal B}$ we have
$\frac{n^2}{2}+2nm+2m^2\le 4m$, a contradiction.

{\bf Case 3} Let $1\in {\cal A}$, $1\notin {\cal B}$. As mentioned
above, the given decomposition induces the following
decompositions of the ideals of the even component:
$$ H(F_n)=\pi_1({\cal A}_0)+\pi_1({\cal B}_0),\eqno (10)$$
$$ H({\cal Q}_m)=\pi_2({\cal A}_0)+\pi_2({\cal B}_0).\eqno (11) $$
If either  $\pi_1({\cal A}_0)$, $\pi_2({\cal B}_0)$ or
$\pi_1({\cal B}_0)$, $\pi_2({\cal A}_0)$ are non-simple
semisimple, then $\dim{\cal A}_0=\dim \pi_1({\cal A}_0)<\dim
H(F_n)$, $\dim {\cal B}_0=\dim\pi_2({\cal B}_0)<\dim H({\cal
Q}_m)$. This implies that $\dim{\cal A}_0+\dim {\cal B}_0<\dim
(H(F_n)\oplus H({\cal Q}_m))$, which is wrong. Likewise we have a
contradiction in the second case. Therefore, there is a simple
algebra in each pair:($\pi_1({\cal A}_0)$, $\pi_2({\cal
B}_0)$),\,($\pi_1({\cal B}_0)$, $\pi_2({\cal A}_0)$). Since 1 is
not an element of $\cal B$, $\cal B$ has a non-zero  two-sided
annihilator, and so does ${\cal B}_0$. It follows that one of
$\pi_1({\cal B}_0)$, $\pi_2({\cal B}_0)$ has a non-zero two-sided
annihilator. Let us assume the first possibility, that is,
$\pi_1({\cal B}_0)$ can be embedded in the simple subalgebra which
also has a non-zero annihilator. Since $H(F_n)$ cannot be written
as the sum of two simple subalgebras, $\pi_1({\cal A}_0)$ should
be non-simple semisimple. This implies that
$$\pi_1({\cal A}_0)\cong \left\{ \begin{array} {c}
                  H(F_k)\oplus H({\cal Q}_l),\,\mbox{или}\\
                  H({\cal R}_k)\oplus H({\cal R}_l)
                   \end{array}\right.\eqno (12)
$$
In other words, we represent $H(F_n)$ as the sum of a non-simple
semisimple subalgebra of form (12) and a subalgebra which has a
non-zero two-sided annihilator.

Let $V$ denote the $n$-column vector space. Then, there exists a
non-zero vector $v\in V$ annihilated by the second subalgebra. By
Lemma 1.2, $\dim\,H(F_n)v=n$. It follows from (10) that
$\dim\pi_1({\cal A}_0)v=n$.

If $\pi_1({\cal A}_0)\cong H(F_k)\oplus H({\cal Q}_l)$, then by
some automorphism of $F_n$ it can be reduced to the following
form:
$$ \left(\begin{array}{cccccc}
         X&  \ldots   &      0 &         0 &  \ldots &       0\\
\vdots    &    \ddots &  \vdots&  \vdots   &  \ddots &   \vdots\\
   0      &  \ldots   &      X &         0 &  \ldots &        0\\
   0      & \ldots    &      0 &         Y &  \ldots &        0\\
\vdots    & \ddots    & \vdots &   \vdots  &  \ddots &       \vdots \\
   0      & \ldots    &      0 &         0 &  \ldots &        Y
   \end{array}\right),
$$
where $X$ is a symmetric matrix of order $k$, $Y$ is a symplectic
matrix of order $2l$. Let $v$ be equal to
$(v_{11},\ldots,v_{1k_1}, v_{21},\ldots,v_{2l_1})^t$ where
$v_{i1}$ is a vector of dimension $k$, $i=1,\ldots,k_1$, $v_{2j}$
is a vector of dimension $2l$, $j=1,\ldots,l_1$. Since
$\pi_1({\cal A}_0)$ contains 1, $kk_1+2ll_1=n$. Then, $\dim\{
Xv_{1i}|\, X\in H(F_k)\}=k$, $\dim\{ Yv_{2j}|\, Y\in H({\cal
Q}_l)\}=2l-1$ (see Lemma 1.4). Therefore, $\dim\pi_1({\cal
A}_0)v=kk_1+(2l-1)l_1<n$, a contradiction. If $\pi_1({\cal
A}_0)\cong H({\cal R}_k)\oplus H({\cal R}_l)$, then by some
automorphism of $F_n$ it can be reduced to
$$
\left(\begin{array}{ccccc}
               X&     0    &      \ldots   &               0 &  0 \\
               0&    X^t   &      \ldots   &               0 &  0\\
         \vdots &   \vdots &      \ddots   &       \vdots    & \vdots\\
               0&      0   &      \ldots   &             Y   &  0\\
               0&      0   &      \ldots   &             0   &  Y^t
           \end{array}\right).\eqno (13)
$$
Hence, $\pi_2({\cal B}_0)$ has a non-trivial two-sided
annihilator, that is, can be embedded in the simple algebra with a
non-zero annihilator. Therefore, $\pi_2({\cal A}_0)$ is non-simple
semisimple because $H({\cal Q}_m)$ cannot be written as the sum of
two simple subalgebras one of which  has a non-zero two-sided
annihilator (see \cite{tm}). As a result, we have the
decomposition of the form: $H({\cal
 Q}_m)=\pi_2({\cal A}_0)+\pi_2({\cal B}_0)$ which in turn induces
 the following
$$ F_{2m}=F_{2m-1}+\langle \pi_2({\cal A}_0)\rangle,$$
in which the first subalgebra clearly has a non-zero two-sided
annihilator, and the second is non-simple semisimple. According to
\cite{BK}, such decomposition cannot exist. The Lemma is proved.
\end{proof}

\begin{lemma}
A superalgebra $\cal J$ of the type $osp(n,m)$ where $n,m>0$
cannot be decomposed into the sum of two proper simple
subsuperalgebras one of which has either the type $K_3$ or $D_t$.
\end{lemma}
\begin{proof}
First we identify $\cal J$ with $osp(n,m)$. Since Kaplansky's
superalgebra $K_3$ has no finite-dimensional associative
specializations, $K_3$ cannot be a subsuperalgebra of a
superalgebra of the type $osp(n,m)$. Assume that $osp(n,m)={\cal
A}+{\cal B}$ where, for example, ${\cal A}\cong D_t$. Then, the
above decomposition induces the following:
$$ H(F_n)=\pi_1(Fe_1\oplus Fe_2)+\pi_1({\cal B}_0),$$
$$ H({\cal Q}_m)=\pi_2(Fe_1\oplus Fe_2)+\pi_2({\cal B}_0).$$

Let us note that  $\dim\pi_1({\cal B}_0)\le \frac{n(n-1)}{2}$ if
it is a simple subalgebra and $\dim\pi_1({\cal
B}_0)\le\frac{n^2-3n+2}{2}+2$ if it is a non-simple semisimple
subalgebra. This implies that $\dim H(F_n)=\frac{n(n+1)}{2}\le
4+\frac{n^2}{2}-\frac{3n}{2}+1$, $2n\le 5$, $n\le 2$. Similarly,
$\dim\pi_2({\cal B}_0)\le 2m^2-5m+3$ if $\pi_2({\cal B}_0)$ is a
simple subalgebra, and $\dim\pi_2({\cal B}_0)\le 2m^2-5m+4$ if
$\pi_2({\cal B}_0)$ is a non-simple semisimple subalgebra. Thus
$\dim H({\cal Q}_m)=2m^2-m\le 2+2m^2-5m+4$, $2m\le 6$,
$m\le\frac{3}{2}$. Therefore, either ${\cal J}\cong osp(1,1)$ or
${\cal J}\cong osp(2,1)$. Since $\dim osp(1,1)=4$, the first case
is not possible. Let $\cal A$ be isomorphic to $D_t$. By the
dimension argument, either ${\cal B}\cong M_{1,1}(F)$ or
$osp(1,1)$.

In turn, $\cal A$ acts completely reducibly in the 4-dimensional
column vector space $V$ because $\dim\, V_0=\dim\,V_1=2$ (see
Theorem 2.3). Moreover, $V=W_1\oplus W_2$, $\dim\, W_1=1$, $\dim
W_2=3$, $W_1$, $W_2$ are invariant subspaces with respect to the
action of $\cal A$. Besides, $\cal A$ acts in $W_1$ trivially and
in $W_2$ irreducibly. It follows that, by some graded automorphism
of $osp(n,m)$, $\cal A$ can be reduced to the following form:
$$    {\cal A}=\left\{  \left(\begin{array}{c|ccc}
                  0&  &\ldots& 0\\
                   \hline
                   &  &        & \\
            \vdots &  &  A'    & \\
                  0&  &        &
                 \end{array}\right)\right\}.
$$
Since $\dim {\cal A}=4$,  ${\cal A}\cong osp(1,1)$. Finally, we
obtain a decomposition of a superalgebra of the type $osp$ into
the sum of two subsuperalgebras of the same type. As proved
before, such decomposition does not exist. The Lemma is proved.
\end{proof}

\section{Decompositions of superalgebras of types $Q(n)$ and
$P(n)$}

First of all we recall the canonical realizations of Jordan
superalgebras of types $P(n)$ and $Q(n)$. A Jordan superalgebra of
the type $Q(n)$ can be represented as the set of all matrices of
order $2n$ which have the following form:
$$     \left\{  \left(\begin{array}{cc}
                A&B\\
                B&A
               \end{array}\right) \right\}
$$
where $A$ and $B$ are any square matrices of order $n$. A
canonical realization of a Jordan algebra of the type $P(n)$
consists of all matrices of the form:
$$        \left\{ \left(\begin{array}{cc}
                A&B\\
                C&A^t
               \end{array}\right)\right\}
$$
where $A$ is any square matrix of order $n$, $B$ is a symmetric
matrix of order $n$, $C$ is a skewsymmetric matrix of order $n$.
Notice here that the even part of $Q(n)$ as well as $P(n)$ is
isomorphic to a simple Jordan algebra of the type $H({\cal R}_n)$.
Besides, $\dim\, Q(n)=\dim\,P(n)=2n^2$.
 Later on only these two properties will be primarily used, hence,
all Lemmas proved in this section are true for Jordan
superalgebras of both types. For definiteness, we consider only
Jordan subalgebras of type $Q(n)$. In several steps, we prove that
no Jordan superalgebras of types $P(n)$ or $Q(n)$ can be
represented as the sum of two proper simple  subsuperalgebras.

\begin{lemma} Let $\cal A$ of the type $osp(p,q)$ be a proper
subsuperalgebra of $Q(n)$. Then $\dim\,{\cal A}\le
\frac{n^2+n}{2}$.
\end{lemma}
\begin{proof}
It follows from the Lemma conditions that ${\cal A}_0\cong
H(F_p)\oplus H({\cal Q}_q)$ is a proper subalgebra of $Q(n)_0$
which is isomorphic to $H({\cal R}_n)$. Therefore, $p+2q\le n$,
$p,q>0$. It is easy to see that the subalgebra takes on its
maximum value when $p+2q=n$. Then $\dim\,{\cal
A}=\frac{p^2+p}{2}+q(2q-1)+2pq$. This implies that $\dim\,{\cal
A}\le\frac{n^2+n}{2}$. The Lemma is proved.
\end{proof}
\begin{lemma}
Let $\cal A$ of the type $M_{k,l}(F)^{(+)}$ where $k,l>0$ be a
proper subsuperalgebras of $Q(n)$. Then $\dim\,{\cal A}\le n^2$.
\end{lemma}
\begin{proof}
Since $\cal A$ is proper, $k+l\le n$, hence $(k+l)^2\le n^2$.
\end{proof}
\begin{lemma} A superalgebra $\cal J$ of either the type $P(n)$ or $Q(n)$,
$n>1$, cannot be represented as the sum of two proper nontrivial
subsuperalgebras one of which has either the type $K_3$ or $D_t$.
\end{lemma}
\begin{proof} For definiteness, we assume that $\cal J$ has the type
$P(n)$. Next, in order to simplify our notation we identify $\cal
J$ with its canonical realization denoted as $P(n)$. By Remark 1
in Section 2 no superalgebra of the type $K_3$ can be a
subsuperalgebra of $P(n)$. Therefore, let $\cal A$ be of type
$D_t$. The given decomposition of $P(n)$ induces that of the form:
$P(n)_0={\cal A}_0+{\cal B}_0$ where ${\cal A}_0=Fe_1\oplus Fe_2$,
$e_1$ and $e_2$ are pairwise orthogonal idempotents. Next we
estimate the dimension of ${\cal B}_0$. If ${\cal B}_0$ is simple,
then $\dim {\cal B}_0\le n^2-2n+1$. If ${\cal B}_0$ is non-simple
semisimple , then $\dim {\cal B}_0\le n^2-2n+2$. As a result,
$\dim P(n)_0=n^2\le 2+n^2-2n+2$, $n\le 2$. The only case which
remains to prove is when $n=2$. By the dimension and rank
arguments, either ${\cal B}\cong M_{1,1}(F)$ or ${\cal B}\cong
osp(1,1)$. In both cases, ${\cal B}_0$ is isomorphic to
$Fe'_1\oplus Fe'_2$ where $e_1'$ and $e_2'$ are pairwise
orthogonal idempotents. Hence $P(2)_0=H({\cal R}_2)=Fe_1\oplus
Fe_2+Fe'_1\oplus Fe'_2$. Both subalgebras contain the identity.
Thus, $\dim H({\cal R}_2)=4>\dim{\cal A}+\dim{\cal B}$. The Lemma
is proved.
\end{proof}

\begin{lemma}
Let $\cal J$ of the type $P(n)$ or $Q(n)$ be represented as the
sum of two proper non-trivial subsuperalgebras $\cal A$ and $\cal
B$ whose even components are semisimple Jordan algebras and one of
them has a non-zero two-sided annihilator. Then ${\cal J}\ne{\cal
A}+{\cal B}$.
\end{lemma}
\begin{proof}
Let ${\cal J}={\cal A}+{\cal B}$, and $\cal A$ have a two-sided
annihilator. Then ${\cal J}_0={\cal A}_0+{\cal B}_0$ where ${\cal
A}_0$, ${\cal B}_0$ are semisimple Jordan subalgebras,
$\mbox{Ann}\,{\cal A}_0\ne\{0\}$. Since ${\cal J}_0\cong H({\cal
R}_n)$, ${\cal J}_0$ can be represented as the set of all matrices
of order $n$, denoted as $F_n^{(+)}$, under the Jordan
multiplication. Obviously, $F_n=\langle {\cal A}_0\rangle +
\langle {\cal B}_0\rangle $ where $\langle {\cal A}_0\rangle $ and
$\langle {\cal B}_0\rangle $ denote associative enveloping
algebras for ${\cal A}_0$ and ${\cal B}_0$, respectively. This
implies that $F_n$ can be written as the sum of two semisimple
subalgebras $\langle {\cal A}_0\rangle $ and  $\langle {\cal
B}_0\rangle $ one of which has a non-zero two-sided annihilator.
This contradicts Proposition 1 in \cite{BK}. The Lemma is proved.
\end{proof}

The following table summarizes all the information obtained above.

\par\medskip
$$\qquad \left| \begin{array}{c|c|c} \hline
   & {\cal A}    & Max\, \, dim\\
  \hline
  1& M_{k,l}(F)^{(+)}    &  n^2\\
  2& osp(p,q)        & \frac{n^2+n}{2}\\
  3& Q(k)     &  2(n-1)^2\\
  4& P(k)     &  2(n-1)^2\\
   \hline
  \end{array}\right|
$$
\par\medskip
\par\medskip

In the second column we list all possible types which
subsuperalgebras of $P(n)$ and $Q(n)$ can have. In the third
column we point out the maximal dimension corresponding to each
subsuperalgera.

\begin{theorem}
Let $\cal J$ have type either $Q(n)$ or $P(n)$, where $n>1$. Then
$\cal J$ cannot be represented as the sum of two proper simple
non-trivial subsuperalgebras $\cal A$ and $\cal B$.
\end{theorem}
\begin{proof} Since the case $P(n)$ is completely similar to the
case $Q(n)$, we give proof only for $Q(n)$. We have the following
cases.

{\bf Case 1.}  $Q(n)={\cal A}+{\cal B}$, ${\cal A}\cong
M_{k,l}(F)^{(+)}$, ${\cal B}\cong M_{s,t}(F)^{(+)}$. Besides, the
dimensions of both subsuperalgebras are not greater than $n^2$.
This implies that the sum in the above decomposition is direct. As
a consequence of this fact, the given decomposition induces the
decomposition of $Q(n)_0$ into the direct sum of two semisimple
subalgebras ${\cal A}_0$ and ${\cal B}_0$ one of which does not
contain the identity of the whole superalgebra or, equivalently,
has a non-trivial two-sided annihilator. By Lemma 4.4, no such
decomposition is possible.

{\bf Case 2.}
 $Q(n)={\cal A}+{\cal B}$, ${\cal A}\cong osp(p,q)$,
           ${\cal B}\cong
           \left\{\begin{array}{c}
            osp(p',q')\\
            M_{k,l}(F)^{(+)}\\
            Q(k)\\
            P(k)
            \end{array}\right..$

Taking into account Lemma 4.1, we can conclude that the
decomposition into the sum of two subsuperalgebras of the type
$osp$ is not possible. Assume that the second decomposition holds
true. According to the above estimates, the dimension of $\cal A$
and $\cal B$ are not greater than $\frac{n^2+n}{2}$ and $n^2$,
respectively. Hence, by the dimension argument, it is not
possible. For the last two cases, the decomposition of the even
part has the form: $H({\cal R}_n)=
 {\cal A}_0+{\cal B}_0$, ${\cal A}_0\cong H(F_p)\oplus H({\cal Q}_q)$,
 ${\cal B}_0\cong H({\cal R}_k)$. If  $1\in {\cal B}_0$, then $k\le\frac{n}{2}$ and $\dim
({\cal A}+{\cal B})\le \frac{n^2+n}{2}+\frac{n^2}{2}<2n^2$.  If
$1\notin H({\cal R}_k)$, then we have a contradiction with Lemma
4.4.

{\bf Case 3.}   $Q(n)={\cal A}+{\cal B}$, ${\cal A}\cong
M_{k,l}(F)^{(+)}$, ${\cal B}\cong \left\{\begin{array}{c}
                                                    Q(m)\\
                                                    P(m)
                                           \end{array}\right..$

This decomposition induces that of the even part: $H({\cal
R}_n)={\cal A}_0+{\cal B}_0$, ${\cal A}_0\cong H({\cal R}_k)\oplus
H({\cal R}_l)$, ${\cal B}_0\cong H({\cal R}_m)$. If  $1\notin
{\cal B}_0$, then again we have a contradiction with Lemma 4.4. If
$1\in {\cal B}_0$, then $k\le\frac{n}{2}$, that is, $\dim {\cal
B}\le\frac{n^2}{2}$. However $\dim\,({\cal A}+{\cal B})\le
n^2+\frac{n^2}{2}<2n^2$, which is wrong.

{\bf Case 4.} Let $Q(n)={\cal A}+{\cal B}$, ${\cal A}\cong P(k)$,
${\cal B}\cong Q(l)$, $k,l<n$. As above, this decomposition
induces the decomposition of the even part $H({\cal R}_n)$ into
the sum of two subalgebras of types $H({\cal R}_k)$ and $H({\cal
R}_l)$. However it follows from the classification of simple
decompositions of simple Jordan algebras that no such
decomposition exists. The Theorem is proved.
\end{proof}

\section{Decompositions of superalgebras of  types $J(V,f)$, $K_3$, $D_t$}

\begin{theorem} Let a superalgebra $\cal J$ have type $J(V,f)$, and
$\cal A$, $\cal B$ be simple non-trivial subsuperalgebras of $\cal
J$. Then ${\cal J}={\cal A}+{\cal B}$ implies that $\cal A$, $\cal
B$ are isomorphic to some superalgebras of non-singular symmetric
bilinear superforms. Moreover, for any decomposition of $V$ into
the sum of nontrivial graded subspaces $V=W_1+W_2$ with
nondegenerate restrictions $f_1$, $f_2$  of $f$ one has
$J(V,f)=J(W_1,f_1)+J(W_2,f_2)$.
\end{theorem}
\begin{proof}
As usual, we identify $\cal J$ with $J(V,f)=(F+V_0)+V_1$ where
$J(V,f)_0=F+V_0$, $J(V,f)_1=V_1$. It follows from the rule of
multiplication, that ${\cal J}_1{\cal J}_1=F1$, where 1 denotes
the identity in $J(V,f)$. In particular, ${\cal A}_1{\cal
A}_1\subseteq F1$ and ${\cal B}_1{\cal B}_1\subseteq F1$. Note
that the idempotents in ${\cal J}_0$ have the form: 1 or
$\frac{1}{2}+v$ where $f(v,v)=\frac{1}{4}$, $v\in V_0$. In
particular, if $v_1$ and $v_2$ are pairwise orthogonal idempotents
in ${\cal J}_0$, $v_1=\frac{1}{2}+v$, $v_2=\frac{1}{2}-v$ where
$v\in V_0$.

Next we consider the following cases.

{\bf (a)} ${\cal A}\cong K_3$ where  $K_3=\langle e,x,y \rangle$,
$[x,y]=e$, $ex=\frac{x}{2}$, $ey=\frac{y}{2}$, $e^2=e$, $e\in
{\cal A}_0$. As mentioned above,  $e=1$ or $e=\frac{1}{2}+v$,
$v\in V_0$. If $e=1$, then, obviously, $ex=x\ne\frac{x}{2}$.
Hence, $e=\frac{1}{2}+v$. Then, on the one hand, $[x,y]=f(x,y)1\ne
0$. Conversely, $[x,y]=\frac{1}{2}+v$ where $v\ne 0$. Hence
$\frac{1}{2}+v=f(x,y)1$. However $f(x,y)1\in F$, $v\in F$, which
is wrong.

{\bf (b)} ${\cal A}\cong D_t$, $t\ne -1,0,1$ where $D_t=\langle
e_1,e_2,x,y \rangle$, $[x,y]=e_1+te_2$, $e_1,e_2$ are pairwise
orthogonal idempotents in  ${\cal A}_0$. Hence there exists $v\in
V_0$, such that  $e_1=\frac{1}{2}+v$, $e_2=\frac{1}{2}-v$. Then
$[x,y]=\frac{1+t}{2}1+(1-t)v$. On the other hand,
$[x,y]=f(x,y)1\ne 0$. However $(1-t)v=0$, that is $t=1$, but
 $D_1\cong J(V',f')$.

{\bf (c)} ${\cal A}\cong osp(n,m)$. Since ${\cal A}_0\le 2$, then
$n=1$, $m=1$. It is easy to check that $osp(1,1)_1\odot
osp(1,1)_1$ cannot be generated by one idempotent.

{\bf (d)} ${\cal A}\cong M_{n,m}(F)$. As in the previous case,
$\mbox{rk}\,{\cal A}_0\le 2$, that is, $n,m=1$. By simple
calculation we can show that $M_{1,1}(F)\odot M_{1,1}(F)$ cannot
be a linear span of one idempotent.

{\bf (e)} ${\cal A}\cong P(n)$  or $ Q(n)$. In this case, the
proof is similar to one in cases  (c) and (d).

It remains to prove the second part of Theorem 5.1, that is, the
existence of the decomposition. A $Z_2$-graded vector space
$V=V_0+V_1$ can be represented as the sum of two $Z_2$-graded
vector subspaces  $W_1$ and $W_2$ in such a way that
$V_0=(W_1)_0+(W_2)_0$ и $V_1=(W_1)_1+(W_2)_1$ and the restriction
of  $f$ to  $(W_1)_0$, $(W_2)_0$, $(W_1)_1$, $(W_2)_1$ are
non-singular. Thus, ${\cal
J}=(F+V_0)+V_1=(F+(W_1)_0+(W_2)_0)+((W_1)_1+(W_2)_1)$ is the sum
of two proper simple subsuperalgebras  of types $J(W_1,f_1)$ and
$J(W_2,f_2)$, respectively.
\end{proof}
\par\medskip
\par\medskip
{\bf Decompositions of $K_3$}

Let $\cal A$  be a subsuperalgebra of $K_3$. Then we have the
following restrictions: $\dim {\cal A}\le 3$ and  $\mbox{rk}\,
{\cal A}_0=1$. Considering all cases one after another, we obtain
that ${\cal A}\cong J(V,f)$ is the only possible case, and $\dim
{\cal A}=2$. Hence $\dim {\cal A}_0=1$ и ${\cal
A}_0=(K_3)_0=\langle e \rangle $, where $e$ is an idempotent. In
other words, $e$ is the identity of $\cal A$. However, if we
consider some element of the form $\alpha x+\beta y$ belonging to
${\cal A}_1$, then $e(\alpha x+\beta y)=\frac{\alpha x+\beta
y}{2}$, that is, $e$ cannot be the identity. This implies that
$K_3$ has no subsuperalgebras of the type $J(V,f)$. As a direct
consequence, we have

\begin{theorem}
A Jordan superalgebra of the type $K_3$ has no decompositions into
the sum of two proper simple non-trivial subsuperalgebras.
\end{theorem}
\par\medskip
\par\medskip
{\bf Decompositions of $D_t$}

Acting in the same manner as in the previous case, we come to the
conclusion that the only possible subsuperalgebras in $D_t$ can be
of types $K_3$ and $J(V,f)$. Let ${\cal A}\cong K_3$ be a
subsuperalgebra  of $D_t$. Then we can choose a basis in $\cal A$
in such a way that ${\cal A}=\langle e', x', y' \rangle $,
$e'x'=\frac{x'}{2}$, $e'y'=\frac{y'}{2}$, $[x',y']=e'$, $e'^2=e'$.
Moreover, ${\cal A}'_0=\langle e'\rangle $, ${\cal A}'_1=\langle
x',y'\rangle$. Therefore, $e$ is an idempotent in $(D_t)_0=\langle
e_1,e_2\rangle$. Hence either $e'=e_i$, $i=1,2$, or $e'=e_1+e_2$.
In the last case, we have $e'(\alpha x+\beta y)=(e_1+e_2)(\alpha
x+\beta y)=(\alpha x+\beta y)\ne \frac{(\alpha x+\beta y)}{2}$.
Hence $e'=e_i$. On the other hand,  $[(\alpha x+\beta y),(\alpha'
x+\beta'y)]=(\alpha\beta'-\beta\alpha')(e_1+te_2)\ne e'=e_i$. This
implies that  $\cal A$ of the type $K_3$ cannot be a
subsuperalgebra $D_t$.

Let ${\cal A}\cong J(V,f)$ be a subsuperalgebra of $D_t$. Then
${\cal A}_0\subseteq (D_t)_0=\langle e_1,e_2\rangle$, ${\cal
A}_1\subseteq (D_t)_1=\langle x,y\rangle$. It is well-known that a
subsuperalgebra of type $J(V,f)$ has the identity $e$, $e\in
\langle e_1,e_2\rangle$, that is, $e=e_1+e_2$. If $\dim{\cal
A}_0>1$, then we always can choose some element of form $(\alpha
e_1+\beta e_2)$ which is linearly independent with $e_1+e_2$, and
$(\alpha e_1+\beta e_2)^2$  is proportionate to $e_1+e_2$. This
implies that $\alpha=\beta$, $\dim{\cal A}_0=1$. However if
$D_t={\cal A}+{\cal B}$ where ${\cal A}\cong J(V_1,f_1)$, ${\cal
B}\cong J(V_2,f_2)$, then  ${\cal A}_0={\cal B}_0=\langle
e_1+e_2\rangle$, that is, $(D_t)_0\ne{\cal A}_0+{\cal B}_0$.

\begin{theorem}
A Jordan superalgebra of the type $D_t$ has no simple
decomposition into the sum of two non-trivial subsuperalgebras.
\end{theorem}

\section{Acknowledgment}

The author uses this opportunity to thank her supervisor Prof.
Bahturin for his helpful cooperation, many useful  ideas and
suggestions.

\end{document}